\documentclass{amsart}
\usepackage{amssymb}

\usepackage{amscd}

\input{tcilatex}
\input tcilatex

\begin{document}
\author{Marek Kara\'{s}}
\title{Tame automorphisms of $\Bbb{C}^{3}$ with multidegree of the form $%
(3,d_{2},d_{3})$}
\date{}
\maketitle

\begin{abstract}
In this note we prove that the sequence $(3,d_{2},d_{3}),$ where $d_{3}\geq
d_{2}\geq 3,$ is the multidegee of some tame automorphism of $\Bbb{C}^{3}\,$%
if and only if $3|d_{2}$ or $d_{3}\in 3\Bbb{N}+d_{2}\Bbb{N}.$
\end{abstract}

\section{Introduction}

By multidegree of a polynomial mapping $F=(F_{1},\ldots ,F_{n}):\Bbb{C}%
^{n}\rightarrow \Bbb{C}^{n},$ denoted $\limfunc{mdeg}F,$ we call the
sequence $(\deg F_{1},\ldots ,\deg F_{n}).$ It seems to be interesting for
which sequences $(d_{1},\ldots ,d_{n})\,$there are automorphisms or tame
automorphisms of $\Bbb{C}^{n}$ with $\limfunc{mdeg}F=(d_{1},\ldots ,d_{n}).$
Let us recall that a tame automorphism is a composition of linear and
triangular automorphisms.

By $\limfunc{Tame}(\Bbb{C}^{n})$ we will denote the group of all tame
automorphisms of $\Bbb{C}^{n}$ and by $\limfunc{mdeg}$ the mapping from the
set of all endomorphisms of $\Bbb{C}^{n}$ into the set $\Bbb{N}^{n}.$

In \cite{Karas} was proven that $(3,4,5),(3,5,7),(4,5,7),(4,5,11)\notin 
\limfunc{mdeg}(\limfunc{Tame}(\Bbb{C}^{3}))$ and that for all $d_{3}\geq
d_{2}\geq 2,\,(2,d_{2},d_{3})\in \limfunc{mdeg}(\limfunc{Tame}(\Bbb{C}%
^{3})). $ Next in \cite{Karas2} it was proven that if $d_{3}\geq
d_{2}>d_{1}>2,$ and $d_{1},d_{2}$ are prime numbers, then $%
(d_{1},d_{2},d_{3})\in \limfunc{mdeg}(\limfunc{Tame}(\Bbb{C}^{3}))$ if and
only if $d_{3}\in d_{1}\Bbb{N}+d_{2}\Bbb{N}.$ In this paper we investigate
the set 
\begin{equation*}
\{(3,d_{2},d_{3})\ |\ 3\leq d_{2}\leq d_{3}\ \}\cap \limfunc{mdeg}(\limfunc{%
Tame}(\Bbb{C}^{3})).
\end{equation*}
Namely we show the following theorem.

\begin{theorem}
\label{main}If $3\leq d_{2}\leq d_{3},$ then $(3,d_{2},d_{3})\in \limfunc{%
mdeg}(\limfunc{Tame}(\Bbb{C}^{3}))$ if and only if $3|d_{2}$ or $d_{3}\in 3%
\Bbb{N}+d_{2}\Bbb{N}.$
\end{theorem}

Since for all permutation $\sigma $ of the set $\{1,2,3\},$ $%
(d_{1},d_{2},d_{3})\in \limfunc{mdeg}(\limfunc{Tame}(\Bbb{C}^{3}))$ if and
only if $(d_{\sigma (1)},d_{\sigma (2)},d_{\sigma (3)})\in \limfunc{mdeg}(%
\limfunc{Tame}(\Bbb{C}^{3})),\,$then the assumption $3\leq d_{2}\leq d_{3}$
is not restrictive.

\section{Some useful results}

For the convenient of the reader we collect in this section all results that
we will need in the proof of Theorem \ref{main}.

\begin{theorem}
\label{tw_sywester}If $a,b$ are positive integers such that $\gcd (a,b)=1,$
then for every integer $k\geq (a-1)(b-1)$ there are $k_{1},k_{2}\in \Bbb{N}$
such that 
\begin{equation*}
k=k_{1}a+k_{2}b.
\end{equation*}
Moreover $(a-1)(b-1)-1\notin a\Bbb{N}+b\Bbb{N}.$
\end{theorem}

In the proof we will, also, use the following proposition.

\begin{proposition}
\textit{(\cite{Karas}, Proposition 2.2) }\label{prop_sum_d_i}If for a
sequence of integers $1\leq d_{1}\leq \ldots \leq d_{n}$ there is $i\in
\{1,\ldots ,n\}$ such that 
\begin{equation*}
d_{i}=\sum_{j=1}^{i-1}k_{j}d_{j}\qquad \text{with }k_{j}\in \Bbb{N},
\end{equation*}
then there exists a tame automorphism $F$ of $\Bbb{C}^{n}$ with $\limfunc{%
mdeg}F=(d_{1},\ldots ,d_{n}).$
\end{proposition}

\begin{definition}
\textit{(\cite{sh umb1}, Definition 1) }\label{def_*-red}A pair $f,g\in
k[X_{1},\ldots ,X_{n}]$ is called *-reduced if\newline
(i) $f,g$ are algebraically independent;\newline
(ii) $\overline{f},\overline{g}$ are algebraically dependent, where $%
\overline{h}$ denotes the highest homogeneous part of $h$;\newline
(iii) $\overline{f}\notin k[\overline{g}]$ and $\overline{g}\notin k[%
\overline{f}].$
\end{definition}

\begin{definition}
\textit{(\cite{sh umb1}, Definition 1) }Let $f,g\in k[X_{1},\ldots ,X_{n}]$
be a *-reduced pair with $\deg f<\deg g.$ Put $p=\frac{\deg f}{\gcd (\deg
f,\deg g)}.$ In this situation the pair $f,g$ is called $p-$reduced pair.
\end{definition}

\begin{theorem}
\textit{(\cite{sh umb1}, Theorem 2)}\label{tw_deg_g_fg} Let $f,g\in
k[X_{1},\ldots ,X_{n}]$ be a $p-$reduced pair, and let $G(x,y)\in k[x,y]$
with $\deg _{y}G(x,y)=pq+r,0\leq r<p.$ Then 
\begin{equation*}
\deg G(f,g)\geq q\left( p\deg g-\deg g-\deg f+\deg [f,g]\right) +r\deg g.
\end{equation*}
\end{theorem}

In the above theorem $[f,g]$ means the Poisson bracket of $f$ and $g,$ but
for us it is only important that 
\begin{equation*}
\deg [f,g]=2+\underset{1\leq i<j\leq n}{\max }\deg \left( \frac{\partial f}{%
\partial x_{i}}\frac{\partial g}{\partial x_{j}}-\frac{\partial f}{\partial
x_{j}}\frac{\partial g}{\partial x_{i}}\right)
\end{equation*}
if $f,g$ are algebraically independent, and $\deg [f,g]=0$ if $f,g$ are
algebraically dependent.

Notice, also, that the estimation from Theorem \ref{tw_deg_g_fg} is true
even if the condition (ii) of Definition \ref{def_*-red} is not satisfied.
Indeed, if $G(x,y)=\sum_{i,j}a_{i,j}x^{i}y^{j},$ then, by the algebraic
independence of $\overline{f}$ and $\overline{g}$ we have: 
\begin{eqnarray*}
\deg G(f,g) &=&\underset{i,j}{\max }\deg (a_{i,j}f^{i}g^{j})\geq \deg
_{y}G(x,y)\cdot \deg g= \\
&=&(qp+r)\deg g\geq q(p\deg g-\deg f-\deg g+\deg [f,g])+r\deg g.
\end{eqnarray*}
The last inequality is a consequence of the fact that $\deg [f,g]\leq \deg
f+\deg g.$

We will also use the following theorem.

\begin{theorem}
\label{tw_type_1-4}\textit{(\cite{sh umb1}, Theorem 3) }Let $%
F=(F_{1},F_{2},F_{3})\,$be a tame automorphism of $\Bbb{C}^{3}.$ If $\deg
F_{1}+\deg F_{2}+\deg F_{3}>3$ (in other words if $F$ is not a linear
automorphism), then $F$ admits either an elementary reduction or a reduction
of types I-IV (see \cite{sh umb1} Definitions 2-4).
\end{theorem}

Let us, also, recall that an automorphism $F=(F_{1},F_{2},F_{3})$ admits an
elementary reduction if there exists a polynomial $g\in \Bbb{C}[x,y]$ and a
permutation $\sigma $ of the set $\{1,2,3\}$ such that $\deg (F_{\sigma
(1)}-g(F_{\sigma (2)},F_{\sigma (3)}))<\deg F_{\sigma (1)}.$

\section{Proof of the theorem}

\begin{proof}
By Corollary \ref{prop_sum_d_i} if $3|d_{2}$ or $d_{3}\in 3\Bbb{N}+d_{2}\Bbb{%
N},$ there exists a tame automorphism $F:\Bbb{C}^{3}\rightarrow \Bbb{C}^{3}$
such that $\limfunc{mdeg}F=(3,d_{2},d_{3}).$ Thus in order to prove Theorem 
\ref{main} it is enough to show that if $3\nmid d_{2}$ and $d_{3}\notin 3%
\Bbb{N}+d_{2}\Bbb{N},$ then there is no tame automorphism of $\Bbb{C}^{3}$
with multidegree $(3,d_{2},d_{3}).$ so from now we will assume that $3\nmid
d_{2}$ and $d_{3}\notin 3\Bbb{N}+d_{2}\Bbb{N}.$

Since $3\nmid d_{2},$ $\gcd (3,d_{2})=1.$ Then by Theorem \ref{tw_sywester},
for all $k\geq (3-1)(d_{2}-1)=2d_{2}-2$ we have $k\in 3\Bbb{N}+d_{2}\Bbb{N}.$
Thus, since $d_{3}\notin 3\Bbb{N}+d_{2}\Bbb{N},$ we have 
\begin{equation}
d_{3}<2d_{2}-2  \label{row_sylwester}
\end{equation}

Assume that $F=(F_{1},F_{2},F_{3})$ is an automorphism of $\Bbb{C}^{3}$ such
that $\limfunc{mdeg}F=(3,d_{2},d_{3}).$ Our aim is to prove that this
hypothetical automorphism can not be a tame automorphism. By Theorem \ref
{tw_type_1-4} it is enough to show that $F$ does not admit neither reduction
of types I-IV (see \cite{sh umb1}, Definitions 2-4) nor elementary reduction.

Assume that $F$ admits a reduction of type I. Then by the definition (see 
\cite{sh umb1}, Definition 2) there is a permutation $\sigma \,$of the set $%
\{1,2,3\}$ and $n\in \Bbb{N}\backslash \{0\}$ such that $\deg F_{\sigma
(1)}=2n,\deg F_{\sigma (2)}=ns,$ where $s\geq 3$ is odd number,\thinspace $%
2n<\deg F_{\sigma (3)}\leq ns.$ Thus we have $2n=d_{2}$ or $2n=d_{3},$ and
then $n\geq 2.$ Since $ns\geq 6>3$ and $2n\geq 4>3,$ then we obtain a
contradiction.

Assume that $F$ admits a reduction of type II. Then by the definition (see 
\cite{sh umb1}, Definition 3) there is a permutation $\sigma \,$of the set $%
\{1,2,3\}$ and $n\in \Bbb{N}\backslash \{0\}$ such that $\deg F_{\sigma
(1)}=2n,\deg F_{\sigma (2)}=3n,$ \thinspace $\frac{3}{2}n<\deg F_{\sigma
(3)}\leq 2n.$ Thus, as before, we have $2n=d_{2}$ or $2n=d_{3},$ and then $%
n\geq 2.$ Since $ns\geq 6>3$ and $\frac{3}{2}n\geq 3,$ then we obtain a
contradiction.

Now assume that $F$ admits a reduction of type III or IV. Then by the
definition (see \cite{sh umb1} Definition 4) there is a permutation $\sigma
\,$of the set $\{1,2,3\}$ and $n\in \Bbb{N}\backslash \{0\}$ such that $\deg
F_{\sigma (1)}=2n,$ and either: 
\begin{equation}
\deg F_{\sigma (2)}=3n,\qquad n<\deg F_{\sigma (3)}\leq \frac{3}{2}n
\label{row_type_III}
\end{equation}
or 
\begin{mathletters}
\begin{equation}
\frac{5}{2}n<\deg F_{\sigma (2)}\leq 3n,\qquad \deg F_{\sigma (3)}=\frac{3}{2%
}n  \label{row_type_IV}
\end{equation}
As before we have $2n=d_{2}$ or $2n=d_{3},$ and $n\geq 2.$ Assume for a
moment that $n>2.$ Then, since $3n,\frac{5}{2}n,\frac{3}{2}n,n+1>3,$ then we
obtain a contradiction. Thus we can assume that $n=2.$ If we assume that (%
\ref{row_type_III}) is hold, then we obtain $d_{2}=4$ and $d_{3}=6.$
Similarly, if we assume that (\ref{row_type_IV}) is hold. This is a
contradiction with $d_{3}\notin 3\Bbb{N}+d_{2}\Bbb{N}.$

Now, assume that $(F_{1},F_{2},F_{3}-g(F_{1},F_{2})),$where $g\in \Bbb{C}%
[x,y],$ is an elementary reduction of $(F_{1},F_{2},F_{3}).$ Hence we have $%
\deg g(F_{1},F_{2})=\deg F_{3}=d_{3}.$ Since $\gcd (3,d_{2})=1,$ then by \ref
{tw_deg_g_fg}, we have 
\end{mathletters}
\begin{equation*}
\deg g(F_{1},F_{2})\geq q(3d_{2}-d_{2}-3+\deg [F_{1},F_{2}])+rd_{2},
\end{equation*}
where $\deg _{y}g(x,y)=3q+r$ with $0\leq r<3.\,$ Since $F_{1},F_{2}$ are
algebraically independent, $\deg [F_{1},F_{2}]\geq 2$ and then $%
3d_{2}-d_{2}-3+\deg [F_{1},F_{2}]\geq 2d_{2}-1.$ Then by (\ref{row_sylwester}%
) follows that $q=0.$ Also by (\ref{row_sylwester}) we must have $r<2.$ Thus 
$g(x,y)=g_{0}(x)+g_{1}(x)y.$ Since $3\Bbb{N\cap (}d_{2}+3\Bbb{N})=\emptyset
, $ then $\deg g(F_{1},F_{2})\in 3\Bbb{N\cup (}d_{2}+3\Bbb{N})\subset 3\Bbb{N%
}+d_{2}\Bbb{N}.$ This is a contradiction.

Now, assume that $(F_{1},F_{2}-g(F_{1},F_{3}),F_{3}),$where $g\in \Bbb{C}%
[x,y],$ is an elementary reduction of $(F_{1},F_{2},F_{3}).$ Therefore we
have $\deg g(F_{1},F_{3})=d_{2}.$ Since $d_{3}\notin 3\Bbb{N}+d_{2}\Bbb{N},$ 
$\gcd (3,d_{3})=1.$ Then by Theorem \ref{tw_deg_g_fg} we have 
\begin{equation*}
\deg g(F_{1},F_{3})\geq q(3d_{3}-d_{3}-3+\deg [F_{1},F_{3}])+rd_{3},
\end{equation*}
where $\deg _{y}g(x,y)=3q+r$ with $0\leq r<3.$ Since $3d_{3}-d_{3}-3+\deg
[F_{1},F_{3}]\geq 2d_{3}-1>p_{2},$ then $q=0.$ Since, also, $d_{3}>d_{2}$
(because $d_{3}\geq d_{2}$ and $d_{3}\notin 3\Bbb{N}+d_{2}\Bbb{N)},$ then $%
r=0.$ Thus $g(x,y)=g(x),$ and $\deg g(F_{1},F_{3})=\deg g(F_{1})\in 3\Bbb{N}%
. $ This is a contradiction with $3\nmid d_{2}.$

Finally, assume that $(F_{1}-g(F_{2},F_{3}),F_{2},F_{3}),$ is an elementary
reduction of $(F_{1},F_{2},F_{3}).$ Thus we have $\deg g(F_{2},F_{3})=3.$
Let 
\begin{equation*}
p=\frac{d_{2}}{\gcd (d_{2},d_{3})}.
\end{equation*}
Since $d_{3}\notin 3\Bbb{N}+d_{2}\Bbb{N},$ $d_{2}\nmid d_{3},$ and then $%
p>1. $ By Theorem \ref{tw_deg_g_fg} we have 
\begin{equation*}
\deg g(F_{2},F_{3})\geq q(pd_{3}-d_{2}-d_{3}+\deg [F_{1},F_{3}])+rd_{3},
\end{equation*}
where $\deg _{y}g(x,y)=qp+r$ with $0\leq r<p.$ Since $d_{3}>3,$ then we have 
$r=0.$ consider the case $p\geq 3.$ In this case $pd_{3}-d_{2}-d_{3}+\deg
[F_{1},F_{3}]\geq d_{3}+\deg [F_{1},F_{3}]>3.$ Thus we must have $q=0.$
Hence $g(x,y)=g(x),$ and $3=\deg g(F_{2},F_{3})=\deg g(F_{2})\in d_{2}\Bbb{N}%
.$ This is a contradiction with $d_{2}\neq 3$ (we have assumption that $%
3\nmid d_{2}$). Consider, now, the case $p=2.$ Since $p=2,$ we have, for
some $n\in \Bbb{N},$ $d_{2}=2n$ and $d_{3}=ns,$ where $s\geq 3$ is odd.
Since, also, $d_{2}>3,$ then $n\geq 2.$ This means that $d_{3}-d_{2}\geq 2,$
and that $2d_{3}-d_{3}-d_{2}+\deg [F_{1},F_{3}]=d_{3}-d_{2}+\deg
[F_{1},F_{3}]\geq 4>3.$ Thus, also in this case we have $q=0.$ As before
this leads to a contradiction
\end{proof}

\vspace{1cm}

\textsc{Marek Kara\'{s}\newline
Instytut Matematyki\newline
Uniwersytetu Jagiello\'{n}skiego\newline
ul. \L ojasiewicza 6}\newline
\textsc{30-348 Krak\'{o}w\newline
Poland\newline
} e-mail: Marek.Karas@im.uj.edu.pl

\end{document}